\documentclass{amsproc}
\usepackage{epsf}
\input{cyracc.def}

\newcommand{\As}{\mathcal{A}s}
\newcommand{\Comm}{\mathcal{C}omm}
\newcommand{\Conf}{\operatorname{Conf}}
\newcommand{\End}[1]{\mathcal{E}nd\,_{#1}}
\newcommand{\id}{\operatorname{id}}
\newcommand{\Hom}{\operatorname{Hom}}
\newcommand{\Map}{\operatorname{Map}}
\newcommand{\MM}{\mathcal{M}}
\newcommand{\nc}{\mathbb{C}}
\newcommand{\nr}{\mathbb{R}}
\newcommand{\nz}{\mathbb{Z}}
\newcommand{\OO}{\mathcal{O}}

\font\cyr=wncyi8

\newtheorem{thm}{Theorem}[section]
\newtheorem{prop}[thm]{Proposition}

\theoremstyle{definition}
\newtheorem{df}[thm]{Definition}

\theoremstyle{remark}
\newtheorem{rem}[thm]{Remark}
\newtheorem*{ack}{Acknowledgment}

\numberwithin{equation}{section}

\begin{document}

\title{The Swiss-Cheese Operad}
\author{Alexander A. Voronov}
\address{Department of Mathematics\\ M.I.T.,
2-246\\ 77 Massachusetts Ave.\\ Cambridge, MA 02139-4307}
%\curraddr{}
\email{voronov@math.mit.edu}
%\urladdr{http://www-math.mit.edu/~voronov/}
\dedicatory{Dedicated to Michael Boardman on the occasion of his sixtieth
birthday.\bigskip\\
\hfill {\cyr \cyracc Vorone gde-to Bog poslal kusochek syra}.
\smallskip\\
\hfill I. A. Krylov
\smallskip}
\date{July 8, 1998}
\thanks{Research supported in part by an AMS Centennial Fellowship.}

\begin{abstract}
We introduce a new operad, which we call the Swiss-cheese operad. It
mixes naturally the little disks and the little intervals operads. The
Swiss-cheese operad is related to the configuration spaces of points
on the upper half-plane and points on the real line, considered by
Kontsevich for the sake of deformation quantization. This relation is
similar to the relation between the little disks operad and the
configuration spaces of points on the plane. The Swiss-cheese operad
may also be regarded as a finite-dimensional model of the moduli space
of genus-zero Riemann surfaces appearing in the open-closed string
theory studied recently by Zwiebach. We describe algebras over the
homology of the Swiss-cheese operad.
\end{abstract}

\maketitle

\section*{Introduction}

This is a short note whose modest purpose is to introduce a new
operad, which we call the Swiss-cheese operad. It mixes naturally the
little disks and the little intervals operads, which are J. Peter
May's \cite{may} modifications of the little cubes operad, the
remarkable discovery of J. Michael Boardman and Rainer M. Vogt
\cite{boardman-vogt}. The Swiss-cheese operad is related to the
configuration spaces of points on the upper half-plane and
points on the real line, considered by Maxim Kontsevich \cite{kon} for
the sake of deformation quantization, in the same way that the little disks
operad is related to the configuration spaces of points on the complex
plane. The Swiss-cheese operad may also be regarded as a
finite-dimensional model of the moduli space of genus-zero Riemann
surfaces appearing in the open-closed string theory studied recently
by Barton Zwiebach
\cite{zwiebach}. Our main theorem describes algebras over the homology
of the Swiss-cheese operad. Such algebraic structure is expected to be
found on the physical state space of open-closed string theory.

\begin{ack}
I would like to thank Takashi Kimura and Jim Stasheff for illuminating
discussions. It would be fair to say that this paper is an offspring
of my collaboration with them. Neither has this paper avoided
Stasheff's famous ``red-ink pen'', which has miraculously turned into
an ``electronic red-ink pen'' at the verge of this millennium. I am
also very grateful to Martin Markl for his valuable remarks on modules
over an operad.
\end{ack}

\section{The Swiss-cheese operad}

We will be using the standard terminology of operad theory, see
\cite{boardman-vogt,gk,kriz-may,loday}. Operads will be considered in
different tensor categories, such as those of manifolds, topological
spaces, (graded) vector spaces, complexes of vector spaces, sets, and
even modules over a given operad, depending on the context. The
Swiss-cheese operad resembles the famous little disks operad. However
the pictures for the Swiss-cheese operad look more like Swiss cheese
than those for the little disks operad. The \emph{little disks operad}
is a collection of manifolds $D(n)$, $n \ge 1$, where each $D(n)$ is
the configuration space of $n$ nonoverlapping little disks inside the
standard unit disk on the plane:
\medskip

\centerline{\epsfxsize=1in \epsfbox{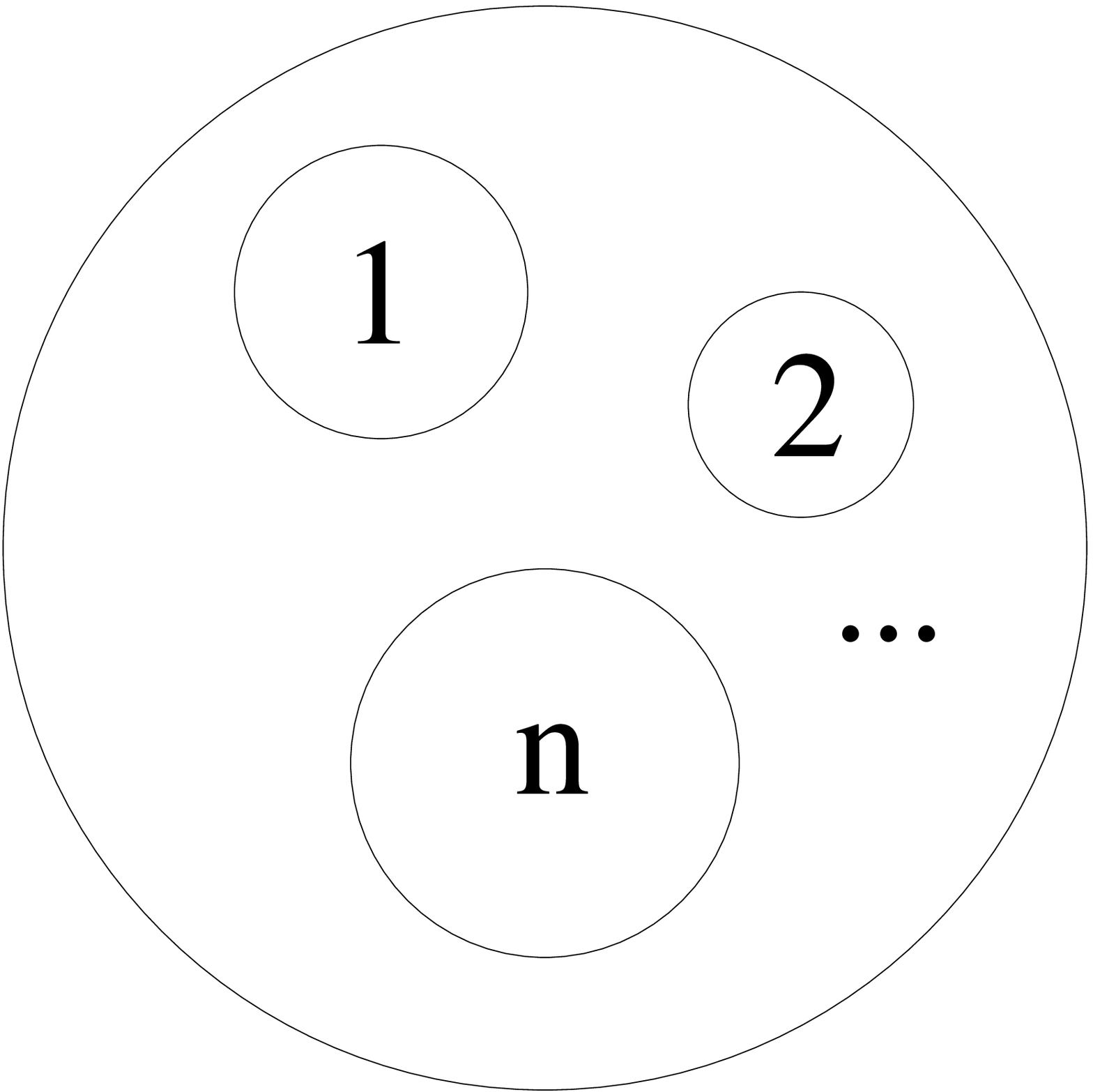}}
\smallskip

The space $D(n)$ is obviously an open set in $\nr^{3n}$ (whence a
manifold structure), each configuration being uniquely determined by
the position of the centers of the disks and their radii. It is assumed
that each little disk is labeled by a number from 1 through $n$, which
defines the action of the permutation group $\Sigma_n$ on $D(n)$. The
operad composition
\[
\gamma: D(k) \times D(n_1) \times \dots \times D(n_k) \to 
D(n_1 + \dots + n_k)
\]
is given by scaling down given configurations in $D(n_1)$, \dots,
$D(n_k)$, gluing them into the $k$ holes in a given configuration in
$D(k)$, and erasing the seams. Thus $D = \{ D(n) \; | \; n \ge 1\}$
becomes an operad of manifolds.

The \emph{Swiss-cheese operad} is a collection of manifolds $S(m,n)$
$m \ge 0$, $n \ge 1$, where $S(m,n)$ is the configuration space of
nonoverlapping disks labeled 1 through $m$ and upper semidisks labeled
1 through $n$ inside the standard unit upper semidisk, so that the
little semidisks are all centered about the diameter of the big
semidisk:
\bigskip

\centerline{\epsfxsize=1.5in \epsfbox{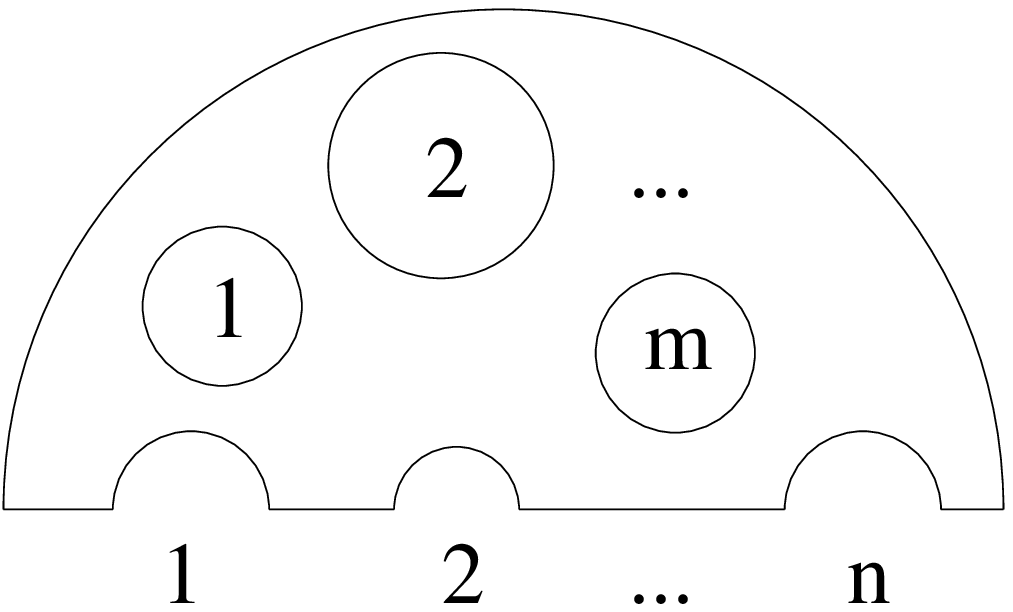}}
\smallskip

$S(m,n)$ is an open subset in $\nr^{3m + 2n}$. There is a natural
action of the group $\Sigma_m \times \Sigma_n$ on $S(m,n)$ and two
types of composition:
\begin{align*}
S(k,l) \times S(m_1,n_1) \times \dots \times S(m_l,n_l)& 
\xrightarrow{\gamma}
S(k + m_1 + \dots + m_l, n_1 + \dots + n_l), \\
S(k,l) \times D(m_1) \times \dots \times D(m_k) & \xrightarrow{\Gamma} 
S(m_1 + \dots + m_k, l).
\end{align*}
The first one is gluing the big semidisks of given configurations in
$S(m_i,n_i)$, $i = 1,\dots, l$, into the $l$ little semidisks of a
given configuration in $S(k,l)$, the second is gluing the big disks of
given configurations in $D(m_i)$, $i = 1, \dots, k$, into the $k$
little disks of a given configuration in $S(k,l)$. This looks like a
new operad-type structure on first sight, but it is nothing but an
\emph{operad over an operad} or a \emph{relative operad}, all meaning
an \emph{operad of modules over a given operad}. The point is that the
collection of all $S(n) = \bigcup_{m \ge 0} S(m,n)$, $n \ge 1$, is an
operad of modules over the little disks operad $D$, each $S(n)$ being
a module over $D$. Recall (\cite{markl,smirnov}) that a (\emph{right})
\emph{module over an operad} $\OO$ is a collection of spaces $\MM(m)$,
$m \ge 0$, on which the action of $\Sigma_m$ is given, as well as maps
\[
\Gamma: \MM(k) \times \OO(m_1) \times \dots \times \OO(m_k) \to 
\MM(m_1 + \dots + m_k)
\]
satisfying virtually the same axioms as those of an operad. Modules
over an operad form a tensor category with respect to the following
tensor product:
\[
(\MM_1 \times \MM_2) (m) = \bigcup_{k + l = m} \Map_{\Sigma_k,\Sigma_l}
(\Sigma_{k+l}, \MM_1(k) \times \MM_2(l)),
\]
where $\Map_{\Sigma_k,\Sigma_l}$ is the set of all $\Sigma_k \times
\Sigma_l$-equivariant maps $\Sigma_k \times \Sigma_l \to \MM_1(k)
\times \MM_2(l)$, $\Sigma_k \times \Sigma_l$ acting on $\Sigma_{k+l}$
by left translation. The tensor product structure is necessary to
speak of operads of $\OO$-modules.

\section{Relative operads and algebras over them}

Another example of a relative operad comes from linear algebra. Let
$V_1$, $V_2$ be two vector spaces. The \emph{endomorphism operad of
the vector space} $V_1$ is the operad $\End{V_1} = \{\Hom
(V_1^{\otimes m}, V_1) \; | \; m \ge 1\}$, the permutation groups
acting by permuting the inputs of multilinear operators and the operad
compositions given by substitution of the outputs into the inputs. The
relative version of the endomorphism operad is the
\emph{relative endomorphism operad} $\End{V_1,V_2} = \{\Hom
(V_1^{\otimes m} \otimes V_2^{\otimes n}, V_2) \; | \; m \ge 0 , n \ge
1\}$. It is obviously an operad of modules over the operad
$\End{V_1}$.

That said, it is easy to define an algebra over a relative operad: it
is analogous to the notion of a group representation. Let $\MM$ be an
operad over an operad $\OO$. An \emph{algebra over the relative
operad} $\MM$ is a pair of vector spaces $V_1$ and $V_2$ and two
mappings
\begin{align*}
\OO & \to \End{V_1},\\
\MM & \to \End{V_1,V_2},
\end{align*}
the former being a morphism of operads, which also makes
$\End{V_1,V_2}$ an $\OO$-module, the latter being a morphism of
relative operads over $\OO$.

Another class of examples of relative operads may be obtained in the
following way. Let $\OO_1$ and $\OO_2$ be two operads such that
$\OO_1$ is an \emph{operad with commutative multiplication}, cf.\
\cite{gv2}, which is nothing but a morphism $\Comm \to \OO_1$ from the
commutative operad $\Comm$ to $\OO_1$. A commutative multiplication
can equivalently be given by an element $M
\in \OO_1(2)$ stable under the action of the symmetric group
$\Sigma_2$ and such that $M \circ M := \gamma(M;M,\id) - \gamma (M;
\id,M) = 0$. (Here if $\OO_1$ is not an operad of abelian groups, the
vanishing of the difference must be understood as the equality of the
two terms in it). We will also assume that the commutative operad has
a \emph{multiplication unit} 1, \emph{i.e}., $\Comm(0)$ is a
point. Therefore, $\OO_1(0)$ will contain a distinguished element as
well. Then the collection $\MM(n) = \OO_1 \times \OO_2(n)= \bigcup_{m
\ge 0} \OO_1(m) \times \OO_2(n)$ for $n \ge 1$ is an operad over
$\OO_1$, the multiplication in $\OO_1$ being used to define a product
\[
\OO_1(k) \times \OO_1(m_1) \times \dots \times \OO_1(m_l) \to 
\OO_1(k + m_1 + \dots + m_l)
\]
needed for defining the operad composition $\gamma$ for $\MM$. The
associativity of the multiplication guarantees the associativity of
the operad composition $\gamma$, and the commutativity implies the
equivariance of $\gamma$ with respect to the action of the symmetric
group.

\begin{prop}
\label{product}
An algebra over the $\OO_1$-operad $\MM = \OO_1 \times \OO_2$ defined
above may be equivalently described as a pair $(V_1, V_2)$, such that
$V_1$ is an algebra over $\OO_1$, and $V_2$ is an algebra over
$\OO_2$, with an action
\begin{align}
\label{ext-prod}
V_1 \otimes V_2 &\to V_2,\\
u \otimes v &\mapsto uv, \notag
\end{align}
defining on $V_2$ an $\OO_2$-algebra structure over the commutative
associative algebra $V_1$, where $V_1$ acquires the structure a
commutative associative algebra with a unit because of the morphism
$\Comm \to \OO_1$. An $\OO_2$-algebra structure on $V_2$ over a
commutative algebra $V_1$ means that
\begin{equation}
\label{parent}
u_1 \dots u_m v
\end{equation}
(where $u_1, \dots, u_m \in V_1$ and $v \in V_2$, $u u'$ denoting the
commutative product in $V_1$, and $uv$ the action of $V_1$ on $V_2$)
is independent of a possible arrangement of parentheses, and
\begin{equation}
\label{linearity}
f(u_1 v_1,\dots, u_n v_n) = u_1 \dots u_{n} f(v_1, \dots, v_n)
\end{equation}
for any $f \in \OO_2(n)$.
\end{prop}

\begin{proof}
1. Given an algebra $(V_1,V_2)$ over the $\OO_1$-operad $\MM = \OO_1
\times \OO_2$, we immediately have an $\OO_1$-algebra structure, which
induces a commutative algebra structure via $\Comm \to \OO_1$, on
$V_1$ and an $\OO_2$-algebra structure on $V_2$. Take the operad units
$\id \in \OO_1(1)$ and $\id \in \OO_2(1)$ and form their cross-product
$(\id,\id) \in \MM$. The algebra structure over the relative operad
$\MM$ will then create a linear map \eqref{ext-prod} corresponding to
$(\id, \id)$. The operad unit axiom for $\OO_1$ and $\OO_2$ implies
that the composite map $(\OO_1(1) \times \OO_2(1) \times \OO_1(m))
\times \OO_2(n) \xrightarrow{\Gamma \times \id} \OO_1(m) \times
\OO_2(1) \times \OO_2(n) \xrightarrow{\gamma} \OO_1(m) \times
\OO_2(n)$ evaluated at $(\id, \id) \in
\OO_1(1) \times \OO_2(1)$ is the identity map. If we track the same
composition for the relative endomorphism operad, we will see that the
operation $V_1^m \times V_2^n \to V_2$ corresponding to $(f_1, f_2)
\in \OO_1(m) \times
\OO_2(n)$ is nothing but 
\begin{equation}
\label{deftn}
f_1(u_1, \dots, u_m) f_2(v_1, \dots, v_n).
\end{equation}

\begin{sloppypar}
To show that an expression \eqref{parent} is independent of an
arrangement of parentheses, it suffices to show that
\[
u_1 (u_2 v) = (u_1 u_2) v.
\]
This follows from the fact that the following operad compositions are
equal: $\gamma ((\id, \id); \linebreak[0] (\id, \id)) \linebreak[1] =
\Gamma ((\id, \id); M)$, by definition of $\gamma$ for the product
relative operad $\OO_1 \times
\OO_2$. The relation \eqref{linearity} is satisfied for the similar
reason: the operad composition $\gamma: \OO_1(0) \times \OO_2(n)
\times (\OO_1(1) \times \OO_2(1)) \times \dots \times (\OO_1(1) \times
\OO_2(1)) \to \OO_1(n) \times \OO_2(n)$ is defined by using the
commutative multiplication $n$ times: $M^{n}: \OO_1(0) \times \OO_1(1)
\times \dots \times \OO_1(1) \to \OO_1(n)$.
\end{sloppypar}

2. Suppose now we have an $\OO_1$-algebra $V_1$, which in particular
means $V_1$ is a commutative associative algebra with a unit, and an
$\OO_2$-algebra $V_2$ over the commutative algebra $V_1$, as defined
in the statement of this proposition. We have to show that these data
define a module $(V_1,V_2)$ over the relative operad $\OO_1 \times
\OO_2$. Define the corresponding morphism $\OO_1(m) \times \OO_2(n)
\to \Hom(V_1^m \otimes V_2^n, V_2)$ by \eqref{deftn}. What remains to
see is that this is indeed a morphism of relative operads. The
symmetric group equivariance follows from that for the structure
operad morphisms $\OO_1 \to \End{V_1}$ and $\OO_2 \to \End{V_2}$. The
compatibility with the actions of $\OO_1$ on $\OO_1 \times \OO_2$ and
$\End{V_1,V_2}$ follows from the fact that $\OO_1 \to \End{V_1}$ is a
morphism of operads. The compatibility with the operad compositions
$\gamma$ for $\OO_1 \times \OO_2$ and $\End{V_1,V_2}$ follows from the
fact that $\OO_2 \to \End{V_2}$ is a morphism of operads and the
linearity \eqref{linearity}.
\end{proof}

\section{Swiss-cheese algebras}

Different operads are in general responsible for different algebraic
structures.  The little disks operad defines the important class of
\emph{Gerstenhaber or G-algebras}, which are defined by two operations, a
(dot) product $ab$ and a bracket $[a,b]$, on a graded vector space
$V$, so that the product defines a graded commutative algebra
structure on $V$ and the bracket a graded Lie algebra structure on
$V[1]$, the desuspension of the graded vector space $V = \bigoplus_n
V^n$: $V[1]^n = V^{n+1}$. The bracket must be a graded derivation of
the product in the following exact sense:
\[
[a , bc] = [a,b] c + (-1)^{(\deg a - 1) \deg b} b [a,c] ,
\]
where $\deg a$ denotes the degree of an element $a \in V$.  In other
words, a G-algebra is a specific graded version of a Poisson
algebra. Here and henceforth, we will consider \emph{G-algebras with a
unit} 1, which is an element in $V$ behaving as a unit with respect to
the dot product and such that $[a,1] = 0$ for all $a \in V$.

The way the little disks operad $D$ has relevance to G-algebras is
through its homology operad $H_\bullet(D)$ and the following theorem.
\begin{thm}[F. Cohen \cite{C1}]
\label{cohen}
The structure of a G-algebra on a $\nz$-grad\-ed vector space is
equivalent to the structure of an algebra over the homology little
disks operad $H_\bullet(D)$.
\end{thm}

\begin{rem}
Here we must take into account the $D(0)$ component of the little
disks operad and a unit in a G-algebra. This is an obvious extension
of Cohen's Theorem.
\end{rem}

What is the analogue of this theorem for the Swiss-cheese operad?
\begin{thm}
An algebra over the homology Swiss-cheese operad $H_\bullet(S)$ is a
pair of graded vector spaces $V_1$ and $V_2$, $V_1$ endowed with the
structure of a G-algebra with a unit and $V_2$ with the structure of a
graded associative algebra over $V_1$, regarded as a graded
commutative algebra with respect to its dot product.
\end{thm}

\begin{df}
A \emph{Swiss-cheese algebra} is a pair $(V_1,V_2)$ of vector spaces
with the structure described in the theorem.
\end{df}

\begin{proof}
An algebra over a relative operad consists of a pair of graded vector
spaces $V_1$ and $V_2$. By Cohen's Theorem \ref{cohen}, $V_1$ is a
G-algebra. Each $D$-module $S(n)$, $n \ge 1$, is homotopy equivalent
to $D \times \Sigma_n$ as a $D$-module. Moreover the whole collection $S =
\{S(n) \; | \; n \ge 1\}$ is homotopy equivalent to the collection
$\{D \times \Sigma_n \; | \; n \ge 1\}$ as a $D$-operad in the homotopy
category of topological spaces --- one just has to choose an arbitrary
point in $D(2)$ as a multiplication, which will be associative and
commutative up to homotopy, and provide the collection of the
symmetric groups $\Sigma_n$'s with the natural operad structure. This
operad structure makes $\{\Sigma_n \; | \; n \ge 1\}$ into the associative
operad $\As$, whose algebras are nothing but (graded) associative
ones. Passing to homology, we get exactly the situation of
Proposition~\ref{product}: $(V_1,V_2)$ is an algebra over the product
operad $H_\bullet (D) \times \Sigma_n$, except that we are now in the
category of graded vector spaces, so that all operads and algebras
will be graded. The proposition then describes this algebra structure
as required by this theorem.
\end{proof}

\section{Homotopy Swiss-cheese algebras}

In this section we will outline a geometric construction of what we
will call the homotopy Swiss-cheese operad, similar to Kontsevich's
construction \cite{kon:feynman} of the $A_\infty$-operad.  We will not
provide complete details whenever the material is a straightforward
generalization of the proofs and constructions of
\cite{ksv1,ksv2,kvz}. We will first construct a relative operad
homotopy equivalent to the Swiss-cheese operad. This new operad will
have a filtration compatible with the relative operad structure and
therefore give rise to an operad of spectral sequences converging to
its homology operad. The spectral sequences will collapse at $E^2$,
making $E^1$ an operad of complexes whose homology is the homology
Swiss-cheese operad. It will be this operad $E^1$ that we will call
the homotopy Swiss-cheese operad.

Consider the configuration space $\Conf_{m,n}$ of $m \ge 0$ labeled
distinct points on the upper half-plane and $n \ge 0$ labeled distinct
points on the real line on the complex plane, assuming $2m + n \ge
2$. Also consider the quotient configuration space $C_{m,n} =
\Conf_{m,n}/G$, where $G = \{ z \mapsto a z + b \;
| \; a, b \in \nr, a > 0\}$ is the group of orientation-preserving
affine transformations of the real line, acting freely because of the
condition $2m + n \ge 2$. The spaces $\Conf_{m,n}$ and $C_{m,n}$ are
homotopy equivalent to the component $S(m,n)$ of the Swiss-cheese
operad. Unlike the spaces $S(m,n)$, the configuration spaces
$\Conf_{m,n}$ or $C_{m,n}$ do not form a relative operad. We will
compactify the spaces $C_{m,n}$ \emph{\`a la} Kontsevich \cite{kon},
and more generally Fulton-MacPherson \cite{fm}, to become smooth
manifolds with corners and form a relative operad, homotopy equivalent
to the Swiss-cheese operad, to which we add the component $S(1,0)$ and
from which throw away the component $S(0,1)$ to satisfy the stability
condition $2m + n \ge 2$.

The compactification $\overline{C}_{m,n}$ will be formed pointwise by
stable configurations. A \emph{stable configuration} is a semistable
complex algebraic curve $X$ of genus zero with $m + n$ smooth labeled
punctures, $m, n \ge 0$, $2m + n \ge 2$, and a special $m+n+1$st
smooth puncture labeled $\infty$, along with the following data:
\begin{enumerate}

\item
The choice of a line through $\infty$ along with the choice of
orientation on the line, which allows us to identify the corresponding
irreducible component $X_\infty$ of $X$ with the complex plane up to a
transformation in the group $G$. No punctures or double points are
allowed strictly below this line on the component $X_\infty$.

\item
The choice of oriented lines passing through selected double points on
selected irreducible components (not more than one line on each
component). No punctures or double points are allowed strictly below
each of these lines on the corresponding irreducible components. If a
double point lies on such a line on one irreducible component, a real
line must pass through this double point on the adjacent component.

\item
The choice of a real tangent direction at the point $\infty$ and to
each irreducible component at each double point not lying on a real
line.

\end{enumerate}
On each irreducible component $X_\alpha$, one can count the number
$m_\alpha$ of punctures and double points not lying on the real line
(if any) and the number $n_\alpha$ of punctures and double points on
the real line (if any). The stability condition means by definition
that, for each component with no real line, $m_\alpha - 1 \ge 2$ and,
for each component with a real line, $2m_\alpha + n_\alpha - 1 \ge
2$. Each component is considered up to the group of its conformal
transformations: $G$ for components with the choice of a real line and
the group $ \{ z \mapsto a z + b \; | \; a \in \nr, a > 0, b \in
\nr^2\}$ for components with no lines.

Naively, one can think of a stable configuration as a semistable
complex algebraic curve which is cut into two pieces by a cross
section through a straight line passing through $\infty$ and a number
of double points on different components, with one of the two pieces
completely forgotten (discarded) and tangent directions chosen at each
double point not lying on the cross section.

As a set, the \emph{compactification} $\overline{C}_{m,n}$ may be
defined as the set of stable configurations. The structure of a smooth
manifold with corners may be obtained by translating the point-set
description into a sequence of real blowups of the space of
configurations where not all points must be distinct. The
compactification we have described here is exactly the same as the one
constructed by Kontsevich in \cite{kon}, which he used to prove his
Formality Conjecture; we are just giving a new interpretation of
points of that compactification.

The simpler similar compactifications $\overline{C}_{m,0}$ in which no
real lines are involved form an operad homotopy equivalent to the
little disks operad. The components $\overline{C}_{0,n}$ form an
operad isomorphic (as an operad of manifolds with corners) to the
Stasheff polyhedra operad, which is homotopy equivalent to the
associative operad $\As$. The whole collection $\overline{C}_{m,n}$
makes up a relative operad over the operad $\overline{C}_{m,0}$. The
two relative operad compositions $\gamma$ and $\Gamma$ are given by
joining stable configurations at punctures to form new double
points. The obtained relative operad is homotopy equivalent to the
Swiss-cheese operad.

Let us filter the spaces $\overline{C}_{m,n}$ by the number of double
points:
\[
\emptyset \subset F^0 \subset \dots \subset F^{2m+n-3} = \overline{C}_{m,n},
\]
where $F^p$ consists of stable configurations that have at least
$2m+n-3-p$ double points, $\dim_\nr F^p = p$. The operad compositions
$\gamma$ and $\Gamma$ respect the filtration, because they create
double points. For each space $\overline{C}_{m,n}$, we get a spectral
sequence. These spectral sequences form an operad of spectral
sequences, see \cite{ksv2}, here a relative one, with the $E^1_{p,q} =
H_{p+q}(F^p,F^{p-1}; \nc) = H^{-q} (F^p \setminus F^{p-1}; \nc)$ terms
forming a relative operad of complexes.  The operad of complexes
$E^1$, by analogy with the corresponding constructions for
$\overline{C}_{0,n}$ and $\overline{C}_{m,0}$, may be called the
\emph{homotopy Swiss-cheese operad}. However, it is not known whether
the spectral sequence collapses at $E^2$ and the homotopy Swiss-cheese
operad is indeed a free resolution of the Swiss-cheese operad. A
\emph{homotopy Swiss-cheese algebra} may be defined as an algebra over
the homotopy Swiss-cheese operad. A homotopy Swiss-cheese algebra
comprises a certain algebraic structure on a pair of complexes $V_1$
and $V_2$ of vector spaces, which includes operations similar to those
of a Swiss-cheese algebra and a hierarchy of higher homotopies for
relations satisfied by those operations in a strict Swiss-cheese
algebra, relations between the relations, etc. Since $V_2$ is an
algebra over the operad $\overline{C}_{0,n}$, it carries the structure
of an $A_\infty$-algebra. Since $V_1$ is an algebra over the operad
$\overline{C}_{m,0}$, it carries the structure of a homotopy
G-algebra, in a sense which is different from the several versions
introduced in \cite{gv1,gj} and used in \cite{gv2} and
\cite{kvz}\footnote{The version of a homotopy G-algebra introduced in
\cite{gj} under the name of a homotopy 2-algebra and utilized in
\cite{gv2,kvz} under the name of $G_\infty$-algebra is not accurate;
the construction that justifies the definition contains an error,
which was pointed out to me by D.~Tamarkin.}. The zeroth row of the
term $E^1$ for the operad $\overline{C}_{m,0}$ is nothing but the
$L_\infty$-operad. The last row is nothing but the $C_\infty$-operad
(see \cite{ksv2}). The only ($q=0$) row in $E^1$ for
$\overline{C}_{0,n}$ is the $A_\infty$-operad. Thus the homotopy
Swiss-cheese algebra structure establishes an interplay between these
three homotopy algebraic structures.

\begin{rem}
A similar compactification was used by Zwiebach in his work
\cite{zwiebach} on open-closed string theory. If one generalizes
our previous work \cite{ksv1,kvz}, one may expect to obtain the
algebraic structure of Ward identities encoded as the structure of a
homotopy Swiss-cheese algebra, see also discussion in Stasheff's
contribution \cite{jim:last} to this volume.
\end{rem}

%\bibliographystyle{amsalpha}
%\bibliographystyle
%{/usr/local/lib/texmf/tex/latex/packages/amslatex/classes/amsalpha}
%\bibliography{swiss}

\begin{thebibliography}{KVZ97}

\bibitem[BV73]{boardman-vogt}
J.~M. Boardman and R.~M. Vogt, \emph{Homotopy invariant algebraic structures on
  topological spaces}, Springer-Verlag, Berlin, 1973, Lectures Notes in
  Mathematics, Vol. 347.

\bibitem[Coh76]{C1}
F.~R. Cohen, \emph{The homology of {$\mathcal{C}_{n+1}$}-spaces, {$n\ge0$}},
  The homology of iterated loop spaces (Berlin), Springer-Verlag, Berlin, 1976,
  Lecture Notes in Mathematics, Vol. 533, pp.~207--351.

\bibitem[FM94]{fm}
W.~Fulton and R.~MacPherson, \emph{A compactification of configuration spaces},
  Ann. of Math. (2) \textbf{139} (1994), no.~1, 183--225.

\bibitem[GJ94]{gj}
E.~Getzler and J.~D.~S. Jones, \emph{Operads, homotopy algebra and iterated
  integrals for double loop spaces}, Preprint, Department of Mathematics, MIT;
  Department of Mathematics Northwestern University, March 1994,
  {\texttt{hep-th/9403055}}.

\bibitem[GK94]{gk}
V.~Ginzburg and M.~Kapranov, \emph{Koszul duality for operads}, Duke Math. J.
  \textbf{76} (1994), no.~1, 203--272, Erratum: Duke Math. J. \textbf{80}
  (1995), no.\ 1, 293.

\bibitem[GV95]{gv2}
M.~Gerstenhaber and A.~A. Voronov, \emph{Homotopy ${G}$-algebras and moduli
  space operad}, Internat. Math. Res. Notices (1995), no.~3, 141--153.

\bibitem[KM95]{kriz-may}
I.~K{\v{r}}{\'\i}{\v{z}} and J.~P. May, \emph{Operads, algebras, modules and
  motives}, Ast\'erisque (1995), no.~233, iv+145pp.

\bibitem[Kon94]{kon:feynman}
M.~Kontsevich, \emph{Feynman diagrams and low-dimensional topology}, First
  European Congress of Mathematics, Vol.\ II (Paris, 1992) (Basel), Progr.
  Math., vol. 120, Birkh\"auser, Basel, 1994, pp.~97--121.

\bibitem[Kon97]{kon}
M.~Kontsevich, \emph{Deformation quantization of {P}oisson manifolds, {I}},
  Preprint, IHES, September 1997, \texttt{math.QA/9709180}.

\bibitem[KSV95]{ksv1}
T.~Kimura, J.~Stasheff, and A.~A. Voronov, \emph{On operad structures of moduli
  spaces and string theory}, Comm. Math. Phys. \textbf{171} (1995), no.~1,
  1--25.

\bibitem[KSV96]{ksv2}
T.~Kimura, J.~Stasheff, and A.~A. Voronov, \emph{Homology of moduli of curves
  and commutative homotopy algebras}, The Gelfand Mathematical Seminars,
  1993--1995 (Boston, MA), Gelfand Math. Sem., Birkh\"auser Boston, Boston, MA,
  1996, pp.~151--170.

\bibitem[KVZ97]{kvz}
T.~Kimura, A.~A. Voronov, and G.~J. Zuckerman, \emph{Homotopy {G}erstenhaber
  algebras and topological field theory}, Operads: Proceedings of Renaissance
  Conferences (Hartford, CT/Luminy, 1995) (Providence, RI), Contemp. Math.,
  vol. 202, Amer. Math. Soc., 1997, pp.~305--333.

\bibitem[Lod96]{loday}
J.-L. Loday, \emph{La renaissance des op\'erades}, Ast\'erisque (1996),
  no.~237, Exp.\ No.\ 792, 3, 47--74, S\'eminaire Bourbaki, Vol.\ 1994/95.

\bibitem[Mar96]{markl}
M.~Markl, \emph{Models for operads}, Comm. Algebra \textbf{24} (1996), no.~4,
  1471--1500.

\bibitem[May77]{may}
J.~P. May, \emph{Infinite loop space theory}, Bull. Amer. Math. Soc.
  \textbf{83} (1977), no.~4, 456--494.

\bibitem[Smi82]{smirnov}
V.~A. Smirnov, \emph{On the cochain complex of topological spaces}, Math. USSR
  Sbornik \textbf{43} (1982), 133--144.

\bibitem[Sta98]{jim:last}
J.~Stasheff, \emph{Grafting {B}oardman's cherry trees to quantum field theory},
  Tech. report, Department of Mathematics, University of Pennsylvania, March
  1998, \texttt{math.AT/9803156}, this volume.

\bibitem[VG95]{gv1}
A.~A. Voronov and M.~Gerstenhaber, \emph{Higher-order operations on the
  {H}ochschild complex}, Funktsional. Anal. i Prilozhen. \textbf{29} (1995),
  no.~1, 1--6.

\bibitem[Zwi97]{zwiebach}
B.~Zwiebach, \emph{Oriented open-closed string theory revisited}, Preprint
  MIT-CTP-2644, HUTP-97/A025, MIT, May 1997, \texttt{hep-th/9705241}.

\end{thebibliography}

\providecommand{\bysame}{\leavevmode\hbox to3em{\hrulefill}\thinspace}

\end{document}